\documentclass[11pt,a4paper,fleqn]{article}
\usepackage{a4wide,amsfonts,amsmath,latexsym,amssymb,euscript,graphicx,units,mathrsfs}

\numberwithin{equation}{section}

\usepackage{graphicx}
\usepackage{color}
\usepackage{amssymb}
\usepackage{amssymb}
\usepackage[T1]{fontenc}
\usepackage{latexsym}
\usepackage{xypic}
\usepackage{eufrak}
\usepackage{euscript}
\usepackage{amsfonts,amsmath}
\usepackage{verbatim}
\usepackage{fancyhdr}
\usepackage[english]{babel}
\usepackage{mathrsfs}
\usepackage{units}

\newtheorem{prop}{Proposition}[section]

\newtheorem{theorem}[prop]{Theorem}
\newtheorem{defi}[prop]{Definition}

\newcommand{\E}{\mathbb{E}}

\renewcommand{\geq}{\geqslant}
\def\leq{\leqslant}

\newcommand{\R}{\mathbb{R}}

\newcommand{\Prob}{\mathbb{P}}

\def\1{{\mathbf{1}}}

\def\1{{\mathbf{1}}}
\def\0.5{{\frac{1}{2}}}

\newcommand{\qed}{\nopagebreak\hspace*{\fill}
{\vrule width6pt height6ptdepth0pt}\par}

\begin{document}

 \pagestyle{plain}

\title{\textbf{Statistical inference for Vasicek-type model\\ 
driven by Hermite processes}}   
\author{Ivan Nourdin, T.T. Diu Tran \\ Universit\'{e} du Luxembourg} 
 \author{Ivan Nourdin\footnote{Universit\'{e} du Luxembourg, Unit\'{e} de Recherche en Math\'{e}matiques. E-mail: ivan.nourdin@uni.lu. IN was partially supported by the Grant F1R-MTH-PUL-15CONF (CONFLUENT) at Luxembourg University} , T. T. Diu Tran\footnote{Universit\'{e} du Luxembourg, Unit\'{e} de Recherche en Math\'{e}matiques. E-mail: ttdiu.tran@gmail.com}}
\maketitle

\begin{abstract}
Let $Z$ denote a Hermite process of order $q \geq 1$ and self-similarity parameter $H \in (\frac{1}{2}, 1)$. This process is $H$-self-similar, has stationary increments and exhibits long-range dependence. When $q=1$, it corresponds to the fractional Brownian motion,  whereas it is not Gaussian as soon as $q\geq 2$. 
In this paper, we deal with a Vasicek-type model driven by $Z$, of the form
$dX_t = a(b - X_t)dt +dZ_t$.
Here, $a > 0$ and $b \in \mathbb{R}$ are considered as unknown drift parameters. We provide estimators for $a$ and $b$ based on continuous-time observations. For all possible values of $H$ and $q$, we prove strong consistency and we 
analyze the asymptotic fluctuations. 
\end{abstract}



\section{Introduction}

Our aim in this paper is to introduce and analyze a {\it non-Gaussian} extension of the fractional model considered in the seminal paper \cite{CR} of Comte and Renault  (see also Chronopoulou and Viens \cite{CV}, as well as the motivations and references therein) and used by these authors to model a situation where, unlike the classical Black-Scholes-Merton model, the volatility exhibits long-memory. More precisely, we deal with the drift parameter estimation problem for a Vasicek-type process $X$,  defined as the unique (pathwise) solution to
\begin{equation}\label{fracV}
X_0=0,\quad dX_t = a(b-X_t)dt + dZ^{q, H}_t, \,\,t\geq 0,
\end{equation}
where $Z^{q,H}$ is a Hermite process of order $q\geq 1$ and Hurst parameter $H\in (\frac12,1)$.
Equivalently, $X$ is the process given explicitly by 
\begin{equation}\label{fracV2}
X_t = b(1-e^{-at}) +  \int_0^t e^{-a(t-s)}dZ^{q, H}_s,
\end{equation}
where the integral with respect to $Z^{q, H}$ must be understood in the Riemann-Stieltjes sense.
In (\ref{fracV}) and (\ref{fracV2}), parameters $a>0$ and $b\in\R$ are considered as  (unknown) real parameters. 

Hermite processes $Z^{q, H}$ of order $q\geq 2$ form a class of genuine non-Gaussian generalizations of the celebrated fractional Brownian motion (fBm), this latter corresponding to the case $q=1$. Like the fBm, they are self-similar, have stationary increments and  exhibit long-range dependence. Their main noticeable difference with respect to fBm is that they are {\it not} Gaussian.
For more details about this family of processes, we refer the reader to Section \ref{sec:hermite}.

As we said, one main practical motivation to study this estimation problem is to provide tools to understand volatility modeling in finance. Indeed, any mean-reverting model in discrete or continuous time can be taken as a model for stochastic volatility. Classical texts can be consulted for this modeling idea; also see the research monograph \cite{FouquePapanicolaou}. Our paper proposes extensions in the type of tail weights for these processes. We acknowledge that there is a gap between the results we cover and their applicability. First, for any practical problems in finance, one must consider discrete-time observations, and one must then choose an observation frequency, keeping in mind that very high frequencies (or infinitely high frequency, which is the case covered in this paper) are known to be inconsistent with models that have no jumps. This problem can be addressed by adding a term to account for microstructure noise in some contexts. Since this falls outside of the scope of this article, we omit any details. Second, for most financial markets, volatility is not observed. One can resort to proxies such as the CBOE's VIX index for volatility on the S\&P500 index. See the paper \cite{ChronopoulouViens} for details about the observation frequency which allows the use of the continuous-time framework in partial observation. In this paper, and for most other authors, such as those for the aforementioned research monograph, these considerations are also out of scope.

Our paper is relevant to the literature on parameter estimation for processes with Gaussian and non-Gaussian long-memory processes, including \cite{BarbozaViens2017,ChronopoulouTudorViens2011,CVT,coeurjolly,OSV,SVSPA2018,HN,bHuNualartZhou,istaslang,NNZ,TudorViens2009}. In the finance context, the highly cited paper \cite{Rosetal} investigates the high-frequency behavior of volatility, drawing on ideas in the paper \cite{Rosenbaum} on long-memory parameter estimation, and before this, the 1997 paper \cite{istaslang}, and the 2001 paper \cite{coeurjolly}. The paper \cite{Rosetal} on rough volatility contends that the short-time behavior indicates that the Hurst parameter $H$ of volatility is less than
$\frac12$. Most other authors point towards $H$ being bigger than $\frac12$, as a model for long memory, which falls within our context, since Hermite processes are limited to this long-memory case. There are many other papers where long-memory noises are used to direct Ornstein-Uhlenbeck and other mean-reverting processes, where estimation is a main motivation. The reader can e.g. consult \cite{OSV,SVSPA2018} and the numerous references therein. These papers are always in the Gaussian context, which the current paper extends. A limited number of attempts have been made to cover non-Gaussian noises, notably \cite{CVT,TTV}. With the exception of the aforementioned papers \cite{ChronopoulouTudorViens2011,CVT,istaslang,TudorViens2009}, all these papers, and the current paper, make no attempt to estimate the Hurst parameter $H$. This is a non-trivial task, which leads authors to consider restrictive self-similar frameworks, and makes quantitative estimator asymptotics difficult to obtain. We leave this question out of the scope of this paper, mentioning only the paper \cite{BarbozaViens2017} which, to the best of our knowledge, is the only paper where $H$ is jointly estimated for the fractional Ornstein-Uhlenbeck process, simultaneously with the process's other parameters, though this paper does not include a rate of convergence in the asymptotics, and covers only the Gaussian case. Thus our paper covers new ground, and uncovers intriguing asymptotic behaviors which are not visible in any of the prior literature.

Let us now describe in more details the results we have obtained.
\begin{defi}
Recall from (\ref{fracV})-(\ref{fracV2}) the definition of the Vasicek-type process $X=(X_t)_{t\geq 0}$ driven by the Hermite process $Z^{q, H}$.  
Assume that $q\geq 1$ and $H\in(\frac12,1)$ are known, whereas $a>0$ and $b\in\R$
are unknown.
Suppose that we continuously observe $X$ over the time interval $[0,T]$, $T>0$.
We define estimators for $a$ and $b$ as follows:
\begin{eqnarray}
\widehat{a}_T&=&\left(\frac{\alpha_T}{H\Gamma(2H)}\right)^{-\frac1{2H}}, \quad\mbox{where $\alpha_T=\frac{1}{T}\int_0^T X_t^2dt - \left(\frac{1}{T}\int_0^T X_tdt\right)^2$},\label{alpha_T}\\
\widehat{b}_{T}&=&\frac{1}{T}\int_0^T X_tdt.\notag
\end{eqnarray}
\end{defi}

In order to describe the asymptotic behavior of $(\widehat{a}_T,\widehat{b}_T)$ when $T\to\infty$, we first need the following proposition, which defines the $\sigma\{Z^{q,H}\}$-measurable random variable $G_\infty$ appearing in (\ref{R4})-(\ref{R}) below.
\begin{prop}\label{finfini}
Assume either ($q=1$ and $H>\frac34$) or $q\geq 2$. Fix $T>0$, and let
$U_T=(U_T(t))_{t\geq 0}$ be the process defined as $U_T(t)=\int_0^t e^{-T(t-u)}dZ^{q, H}_u$.
Eventually, define the random variable $G_T$ by
\[
G_T = T^{\frac{2}{q}(1-H)+2H}\int_0^1 \big(U_T(t)^2 - \E[U_T(t)^2])dt.
\]
Then $G_T$ converges in $L^2(\Omega)$ to a limit written $G_\infty$.
Moreover, $G_\infty/B_{H,q}$ is distributed according to the Rosenblatt distribution of parameter $1-\frac{2}{q}(1-H)$, where
\begin{equation}\label{beta}
B_{H,q}=\frac{H(2H-1)}{\sqrt{(H_0-\frac{1}{2})(4H_0-3)}}\,
\times \frac{\Gamma(2H+\frac2q(1-H))}{2H+\frac2q(1-H)-1},\quad\mbox{with }
H_0=1-\frac{1-H}{q}.
\end{equation}
(The definition of the Rosenblatt random variable is recalled in Definition \ref{rosenblatt}.)
\end{prop}

We can now describe the asymptotic behavior of $(\widehat{a}_T,\widehat{b}_T)$ as
$T\to\infty$.
In the limits (\ref{R4}) and (\ref{R}) below, note that the two components are (well-defined and) correlated, because $G_\infty$ is $\sigma\{Z^{q,H}\}$-measurable by construction. 
\begin{theorem}\label{main}
Let $X=(X_t)_{t\geq 0}$ be given by (\ref{fracV})-(\ref{fracV2}), where $Z^{q, H}=(Z^{q, H}_t)_{t\geq 0}$ is a Hermite process of order $q\geq 1$ and parameter $H\in(\frac12,1)$, and where $a>0$ and $b\in\R$
are (unknown) real parameters.
The following convergences take place as $T\to\infty$.
\begin{enumerate}
\item{} {\rm [Consistency]} 
$
(\widehat{a}_T,\widehat{b}_{T})\overset{\rm a.s.}{\to} (a,b).
$
\item{} {\rm [Fluctuations]}
They depend on the values of $q$ and $H$.
\begin{itemize}
\item {\rm (Case $q=1$ and $H<\frac34$)} 
\begin{eqnarray}
\left(\sqrt{T}\{\widehat{a}_T-a\},
T^{1-H}\{\widehat{b}_{T}-b\}
\right)&\overset{\rm law}{\to}&\left(-\frac{a^{1+4H}\sigma_H}{2H^2\Gamma(2H)}\,N,
\frac{1}{a}N'\right),\label{G1}
\end{eqnarray}
where $N,N'\sim \mathcal{N}(0,1)$ are independent and $\sigma_H$ is given by 
\begin{equation}\label{sigma}
\sigma_H=\frac{2H-1}{H\Gamma(2H)^2}\,\sqrt{\int_\R\left(
\int_{\R_+^2}e^{-(u+v)}|u-v-x|^{2H-2}dudv
\right)^2dx}.
\end{equation}
\item {\rm (Case $q=1$ and $H=\frac34$)} 
\begin{equation}\label{34}
\left(\sqrt{\frac{T}{\log T}}\{\widehat{a}_T-a\},T^{\frac14}\big\{\widehat{b}_T - b\} \right)
\to
\left(\frac{3}4\sqrt{\frac{a}\pi}\, N,\frac{1}{a}N'\right),
\end{equation}
where $N,N'\sim \mathcal{N}(0,1)$ are independent.
\item {\rm (Case $q=1$ and $H>\frac34$)} 
\begin{eqnarray}
\left(T^{2(1-H)}\{\widehat{a}_T-a\},T^{1-H}\big\{\widehat{b}_T - b\} \right)\overset{\rm law}{\to}
\left(-\frac{a^{2H-1}}{2H^2\Gamma(2H)}\Big(G_\infty-(B^H_1)^2\Big),\frac{1}{a}B^H_1\right),\notag\\
\label{R4}
\end{eqnarray}
where $B^H=Z^{1,H}$ is the fractional Brownian motion and
where the $\sigma\{B^H\}$-measurable random variable $G_\infty$
is defined in Proposition \ref{finfini}.

\item {\rm (Case $q\geq 2$ and any $H$)} 
\begin{eqnarray}
\left(T^{\frac2q(1-H)}\{\widehat{a}_T-a\},T^{1-H}\big\{\widehat{b}_T - b\} \right)
\overset{\rm law}{\to}
\left(-\frac{a^{1-\frac2q(1-H)}}{2H^2\Gamma(2H)}\,G_\infty,\frac{1}{a}Z^{q, H}_1\right),
\label{R}
\end{eqnarray}
where the $\sigma\{Z^{q,H}\}$-measurable random variable $G_\infty$
is defined in Proposition \ref{finfini}.
\end{itemize}
\end{enumerate}
\end{theorem}

As we see from our Theorem \ref{main}, strong consistency for $\widehat{a}_T$ and $\widehat{b}_T$ always holds, {\it irrespective} of the values of $q$ (and $H$). That is, when one is only interested in the first order approximation for $a$ and $b$, Vasicek-type model (\ref{fracV})-(\ref{fracV2}) displays a kind of universality with respect to the order $q$ of the underlying Hermite process.
But, as point 2 shows, the situation becomes different when one looks at the fluctuations, that is, when one seeks to construct asymptotic confidence intervals: they heavily depend on $q$ (and $H$).
Furthermore, we highlight the dependence of two components in the limit (\ref{R4}) and (\ref{R}), which is very different of
 the case $q=1$ and $H\leq \frac34$.

The rest of the paper is structured as follows. Section 2 presents some
basic results about multiple Wiener-It\^o integrals and Hermite processes, as well as some other facts which are used throughout the paper. The proof of Proposition 1.2 is then given in Section 3. 
Section 4 is devoted to the proof of 
the consistency part of Theorem \ref{main}, whereas
the fluctuations are analyzed in Section 5. 

\section{Preliminaries}

\subsection{Multiple Wiener-It\^{o} integrals}

Let $B=\big\{B(h),\,h\in L^2(\R)\big\}$ be a Brownian field defined on a probability space $(\Omega,\mathcal{F},\mathbb{P})$, that is, a centered Gaussian family
satisfying $\E[B(h)B(g)]=\langle h,g\rangle_{L^2(\R)}$ for any $h,g\in L^2(\R)$.

For every $q\geq 1$, the $q$th Wiener chaos $\mathcal{H}_q$ is defined as the closed linear subspace of $L^2(\Omega)$ generated by the family of random variables $\{H_q(B(h)),\,h\in L^2(\R),\,\|h\|_{L^2(\R)}=1\}$, where $H_q$ is the $q$th Hermite polynomial ($H_1(x)=x$, $H_2(x)=x^2-1$, $H_3(x)=x^3-3x$, and so on).

The mapping $I^B_q(h^{\otimes q})=H_q(B(h))$ can be extended to a linear isometry between $L^2_s(\R^q)$ (= the space of symmetric square integrable functions of $\R^q$, equipped with the modified
norm $\sqrt{q!}\|\cdot\|_{L^2(\R^q)}$) and
the $q$th Wiener chaos  $\mathcal{H}_q$. When $f\in L^2_s(\R^q)$, the random variable $I^B_q(f)$ is called the {\it multiple Wiener-It\^o integral of $f$ of order $q$}; equivalently, 
one may write
\begin{equation}\label{multipleiqbf}
I_q^B(f) = \int_{\mathbb{R}^q} f(\xi_1, \ldots, \xi_q) dB_{\xi_1}\ldots dB_{\xi_q}.
\end{equation}

Multiple Wiener-It\^o integrals enjoy many nice properties. We refer to \cite{bNourdinPeccatibook} or \cite{bNualartbook} for a comprehensive list of them. Here, we only recall the orthogonality relationship, the isometry formula and the hypercontractivity property. 

First, the {\it orthogonality relationship} (when $p\neq q$) or {\it isometry formula} (when $p=q$) states that, if $f \in L^2_s(\mathbb{R}^p)$ and $g \in L^2_s(\mathbb{R}^q)$ with $p,q\geq 1$, then
\begin{equation}\label{isom**}
\E[I_p^B(f)I_q^B(g)] =
\begin{cases}
 & p! \big\langle f, g \big\rangle_{L^2(\mathbb{R}^p)} \qquad\text{if } p=q\\
& 0 \qquad\qquad\qquad\quad \text{ if } p \ne q.
\end{cases}
\end{equation}

Second, the {\it hypercontractivity property} reads as follows: for any $q\geq 1$, any $k\in[2,\infty)$ and any $f\in L^2_s(\R^q)$,
 \begin{equation}\label{eq:hypercontractivity1}
\E[|I_q^B(f)|^k]^{1/k} \leq (k-1)^{q/2}\E[|I_q^B(f)|^2]^{1/2}.
 \end{equation}
As a consequence, for any $q\geq 1$ and any $k\in[2,\infty)$, there exists a constant $C_{k,q}>0$ such that, for any $F \in \oplus_{l=1}^q \mathcal{H}_l$, we have
 \begin{equation}\label{eq:hypercontractivity2}
\E[|F|^k]^{1/k} \leq C_{k, q}\,\sqrt{\E[F^2]}.
 \end{equation}

\subsection{Hermite processes}\label{sec:hermite}

We now give the definition and present some basic properties of Hermite processes. We refer the reader to the recent book \cite{bTudorbook} for any missing proof and/or any unexplained notion.

\begin{defi}
\textit{The Hermite process} $(Z_t^{q, H})_{t \geq 0}$ of order $q \geq 1$ and self-similarity parameter $H \in (\frac{1}{2}, 1)$ is defined as
\begin{equation}\label{bH1}
Z^{q, H}_t= c(H, q) \int_{\mathbb{R}^q} \bigg( \int_0^t \prod_{j=1}^q(s- \xi_j)_+^{H_0 - \frac{3}{2}}ds\bigg) dB_{\xi_1}\ldots dB_{\xi_q},
\end{equation}
where 
\begin{equation}\label{bH2}
c(H, q) = \sqrt{\frac{H(2H - 1)}{q! \beta^q(H_0 - \frac{1}{2}, 2-2H_0)}} \quad \text{and} \quad H_0 = 1+\frac{H-1}{q} \in \left(1-\frac{1}{2q}, 1\right).
\end{equation}
(The integral (\ref{bH1}) is a multiple Wiener-It\^{o} integral of order $q$ of the form (\ref{multipleiqbf}).)
\end{defi}

The positive constant $c(H, q)$ in (\ref{bH2}) has been chosen to ensure that $\E[(Z_1^{q, H})^2] = 1$.

\begin{defi}
A random variable with the same law as $Z^{q, H}_1$ is called a \textit{Hermite random variable} of order $q$ and parameter $H$.
\end{defi}

Hermite process of order $q=1$ is nothing but the fractional Brownian motion. 
It is the only Hermite process to be Gaussian (and that one could have defined for $H\leq\frac12$ as well). For other reasons, the value $q=2$ is also special. Indeed, Hermite process of order 2 has attracted a lot of interest in the recent past (see \cite{taqqu} for a nice introduction and an overview), and has won its own name: it  is called \textit{the Rosenblatt process}.

\begin{defi}\label{rosenblatt}
A random variable with the same law as $Z^{2, H}_1$ is called a \textit{Rosenblatt random variable}.
\end{defi}

Except for Gaussianity, Hermite processes of order $q \geq 2$ share many properties with the fractional Brownian motion (corresponding to $q=1$). We list some of them in the next statement.

\begin{prop}
The Hermite process $Z^{q,H}$ of order $q\geq 1$ and Hurst parameter $H\in (\frac12,1)$ enjoys the following properties.
\begin{itemize} 
\item{} {\rm [Self-similarity]} For all  $c > 0, (Z^{q, H}_{ct})_{t \geq 0} \overset{law}{=}  (c^H Z^{q, H}_t)_{t \geq 0}$. 
\item {}{\rm [Stationarity of increments]} For any $h >0$, $(Z^{q, H}_{t+h} - Z^{q, H}_h)_{t \geq 0} \overset{law}{=} (Z^{q, H}_t)_{t \geq 0}$.
\item {}{\rm [Covariance function]} For all $s, t \geq 0$, $\E[Z_t^{q, H}Z_s^{q,H}]= \frac{1}{2}(t^{2H} + s^{2H} - |t-s|^{2H})$.
\item  {}{\rm [Long-range dependence]} $\sum_{n=0}^\infty |\E[Z_1^{q, H}(Z_{n+1}^{q, H} - Z_n^{q, H})]| = \infty$.
\item {}{\rm [H\"{o}lder continuity]} For any $\zeta \in (0, H)$ and any compact interval $[0,T]\subset\R_+$, $(Z^{q, H}_t)_{t\in [0,T]}$ admits a version with H\"{o}lder continuous sample paths of order $\zeta$.
\item {}{\rm [Finite moments]} For every $p \geq 1$, 
there exists a constant $C_{p,q}>0$ such that
$\E[|Z^{q, H}_t|^p] \leq C_{p, q} t^{pH}$ for all $t\geq 0$.
\end{itemize}
\end{prop}

\subsection{Wiener integral with respect to Hermite processes}

The Wiener integral of a deterministic function $f$ with respect to a Hermite process $Z^{q, H}$, which we denote by $\int_{\mathbb{R}}f(u)dZ^{q, H}_u$, has been constructed by Maejima and Tudor in \cite{bMaejimaTudor}. 

Below is a very short summary of what will is needed in the paper about those integrals.
The stochastic integral $\int_{\mathbb{R}}f(u)dZ^{q, H}_u$  is  well-defined for any $f$ belonging to the space $|\mathcal{H}|$
of functions $f: \mathbb{R} \to \mathbb{R}$ such that 
\[
\int_{\mathbb{R}}\int_{\mathbb{R}} |f(u)f(v)||u-v|^{2H-2}dudv < \infty.
\]
We then have, for any  $f,g\in|\mathcal{H}|$, that
\begin{equation}\label{bisometry}
\E\bigg[\int_{\mathbb{R}}f(u)dZ^{q, H}_u\int_{\mathbb{R}}g(v)dZ^{q, H}_u\bigg] = H(2H-1)\int_{\mathbb{R}}\int_{\mathbb{R}}f(u)g(v)|u - v|^{2H-2}dudv.
\end{equation}
Another important and useful property is that, whenever $f \in |\mathcal{H}|$, the stochastic integral $\int_{\mathbb{R}}f(u)dZ^{q, H}_u$ admits the following representation as a multiple Wiener-It\^{o} integral of the form (\ref{multipleiqbf}):
\begin{equation}\label{beq:13}
\int_{\mathbb{R}}f(u)dZ^{q, H}_u = c(H, q) \int_{\mathbb{R}^q}\bigg(\int_{\mathbb{R}}f(u) \prod_{j=1}^q(u- \xi_j)_+^{H_0 - \frac{3}{2}}du\bigg)dB_{\xi_1}\ldots dB_{\xi_q},
\end{equation}
with $c(H, q)$ and $H_0$ given in (\ref{bH2}). 

\subsection{Existing limit theorems}

To the best of our knowledge, only a few limit theorems have been already obtained in the litterature for quadratic functionals of the Hermite process, see \cite{ChronopoulouTudorViens2011,Diu, TudorViens2009}. Here, we only recall the following result from 
\cite{Diu}, 
which will be useful to study the fluctuations of $(\widehat{a}_T,\widehat{b}_T)$ in Theorem \ref{main}. 

\begin{prop}
Assume either $q\geq 2$ or ($q=1$ and $H>\frac34$), and
let $Y$ be 
defined as
\begin{equation}\label{beq:X}
Y_t = \int_0^t e^{-a(t-u)}dZ^{q, H}_u, \qquad t \geq 0.
\end{equation}
Then, 
as $T\to\infty$,
\begin{equation}\label{diu}
T^{\frac{2}{q}(1-H)-1}\int_0^T \big(Y_t^2-\E[Y_t^2]\big)dt
\overset{\rm law}{\to} B_{H,q}\, a^{-2H-\frac2q(1-H)}\times  R^{H'}_1,
\end{equation}
where $R_1^{H'}$ is distributed according to a Rosenblatt random variable of parameter $H'=1-\frac{2}{q}(1-H)$ and $B_{H,q}$ is given by (\ref{beta}).
\end{prop}

Along the proof of Theorem \ref{main}, we will also make use of another result for the Gaussian case ($q=1$), which we take from \cite{NNZ}.
\begin{prop}
Let $Y$ be given by (\ref{beq:X}), with $q=1$ and $H\in\big(\frac12,\frac34\big)$. Then, 
as $T\to\infty$,
\begin{equation}\label{ivan}
T^{-\frac12}\int_0^T \big(Y_t^2-\E[Y_t^2]\big)dt
\overset{\rm law}{\to} a^{2H}\sigma_H\,N,
\end{equation}
where $\sigma_H$ is given by (\ref{sigma}) and $N\sim \mathcal{N}(0,1)$.
\end{prop}

Since $T^{-\frac12}\int_0^T \big(Y_t^2-\E[Y_t^2]\big)dt$ (resp. $T^{-H}B^H_T$) belongs to the second (resp. first) Wiener chaos, 
we deduce from (\ref{ivan}) and the seminal Peccati-Tudor criterion (see, e.g., \cite[Theorem 6.2.3]{bNourdinPeccatibook}) that
\begin{equation}\label{ivanbis}
\left(T^{-\frac12}\int_0^T \big(Y_t^2-\E[Y_t^2]\big)dt, T^{-H}B^H_T\right)
\overset{\rm law}{\to} (a^{2H}\sigma_H\,N, N'),
\end{equation}
where $N,N'\sim N(0,1)$ are independent.

Finally, in the critical case $q=1$ and $H=\frac34$, we will need the following result,
 established in \cite[Theorem 5.4]{bHuNualartZhou}.

\begin{prop}
Let $Y$ be given by (\ref{beq:X}), with $q=1$ and $H=\frac34$. Then, 
as $T\to\infty$,
\begin{equation}\label{nunu}
(T\log T)^{-\frac12}\int_0^T \big(Y_t^2-\E[Y_t^2]\big)dt
\overset{\rm law}{\to} \frac{27}{64a^2}\,N,
\end{equation}
where  $N\sim N(0,1)$.
\end{prop} 

Similarly to (\ref{ivanbis}) and for exactly the same reason, we actually have 
\begin{equation}\label{nunubis}
\left((T\log T)^{-\frac12}\int_0^T \big(Y_t^2-\E[Y_t^2]\big)dt
, T^{-\frac34}B^{\frac34}_T\right)
\overset{\rm law}{\to} (\frac{27}{64a^2}\,N, N'),
\end{equation}
where $N,N'\sim N(0,1)$ are independent.

\subsection{A few other useful facts}
In this section, we let $X$ be given by (\ref{fracV2}), with $a>0$, $b\in\R$ and $Z^{q, H}$ a Hermite process of order $q\geq 1$ and Hurst parameter $H\in(\frac12,1)$.
We can write
\begin{equation}\label{decompoX}
X_t=h(t)+Y_t,
\quad\mbox{
 where
$h(t)=b(1-e^{-at})$ and $Y_t$ is given by (\ref{beq:X}).}
\end{equation}
The following limit, obtained as a consequence of the isometry property (\ref{bisometry}),  will be used many times throughout the sequel:
\begin{eqnarray}
\E[Y_T^2]&=&H(2H-1)\int_{[0,T]^2} e^{-a(T-u)}e^{-a(T-v)}|u-v|^{2H-2}dudv\notag\\
&=&H(2H-1)\int_{[0,T]^2} e^{-a\,u}e^{-a\,v}|u-v|^{2H-2}dudv\notag\\
&\to&H(2H-1)\int_{[0,\infty)^2} e^{-a\,u}e^{-a\,v}|u-v|^{2H-2}dudv\notag\\
&&=a^{-2H}H\Gamma(2H)<\infty. 
\label{yt2}
\end{eqnarray}
Identity (\ref{yt2}) comes from
\begin{eqnarray}\label{eq:sao}
&&(2H-1)\int_{[0,\infty)^2}e^{-a(t+s)}|t-s|^{2H-2}dsdt \nonumber \\
&=&
a^{-2H}(2H-1)\int_{[0,\infty)^2}e^{-(t+s)}|t-s|^{2H-2}dsdt
= a^{-2H}\Gamma(2H) ,
\end{eqnarray}
see, e.g., Lemma 5.1 in Hu-Nualart \cite{HN} for the second equality.
In particular, we note that 
\begin{equation}\label{O1}
\E[Y_T^2]=O(1)\quad\mbox{ as $T\to\infty$}.
\end{equation}
Another simple but important fact that will be used is the following identity:
\begin{equation}\label{ivan2}
\int_0^T Y_tdt = \frac{1}{a}(Z^{q, H}_T-Y_T),
\end{equation}
which holds true since
\begin{eqnarray*}
\int_0^T Y_tdt&=&\int_0^T \left( \int_0^t e^{-a (t-u)}dZ^{q, H}_u\right)dt=\int_0^T \left(\int_u^T e^{-a(t-u)}dt\right)dZ^{q, H}_u=\frac{1}{a}(Z^{q, H}_T-Y_T).
\end{eqnarray*}

\section{Proof of Proposition \ref{finfini}}

We are now ready to prove Proposition \ref{finfini}.

We start by showing that $G_T$ converges well in $L^2(\Omega)$. In order to do so, we will check that the Cauchy
criterion is satisfied.
According to (\ref{beq:13}), we can write $U_T(t)=c(H,q)I_q(g_T(t,\cdot))$, where
\[
g_T(t,\xi_1,\ldots,\xi_q)=\int_0^t e^{-T(t-v)}\prod_{j=1}^q (v-\xi_j)_+^{H_0-\frac32}dv. 
\]
As a result,  we can write, thanks to \cite[identity (3.25)]{bNourdinRosinski},
\begin{eqnarray*}
&&{\rm Cov}(U_S(s)^2,U_T(t)^2) \\
&=& c(H,q)^4
\sum_{r=1}^q \binom{q}{r}^2\bigg\{
q!^2\|g_S(s,\cdot)\otimes_r g_T(t,\cdot)\|^2
+r!^2(2q-2r)!\|g_S(s,\cdot)\widetilde{\otimes}_r g_T(t,\cdot)\|^2\bigg\},
\end{eqnarray*}
implying in turn that
\begin{eqnarray*}
&&\E[G_TG_S]\\
&=&(ST)^{\frac{2}{q}(1-H)+2H}
\int_{[0,1]^2}{\rm Cov}(U_S(s)^2,U_T(t)^2)dsdt\\
&=&c(H,q)^4(ST)^{\frac{2}{q}(1-H)+2H}
\sum_{r=1}^q \binom{q}{r}^2
q!^2 \int_{[0,1]^2}\|g_S(s,\cdot)\otimes_r g_T(t,\cdot)\|^2dsdt\\
&&+c(H,q)^4(ST)^{\frac{2}{q}(1-H)+2H}
\sum_{r=1}^q \binom{q}{r}^2 r!^2(2q-2r)! \int_{[0,1]^2}\|g_S(s,\cdot)\widetilde{\otimes}_r g_T(t,\cdot)\|^2dsdt.
\end{eqnarray*}
To check the Cauchy criterion for $G_T$, we are thus left to show the existence, 
for any $r\in\{1,\ldots,q\}$, of
\begin{eqnarray}
&&\lim_{S,T\to\infty}(ST)^{\frac{2}{q}(1-H)+2H}
 \int_{[0,1]^2}\|g_S(s,\cdot)\otimes_r g_T(t,\cdot)\|^2dsdt \label{limit1}\\
 &\mbox{and}&\lim_{S,T\to\infty}(ST)^{\frac{2}{q}(1-H)+2H}
 \int_{[0,1]^2}\|g_S(s,\cdot)\widetilde{\otimes}_r g_T(t,\cdot)\|^2dsdt. \label{limit2}
\end{eqnarray}
Using that $\int_{\mathbb{R}} (u-x)_+^{H_0-\frac32}(v-x)_+^{H_0-\frac32}du =c_H\,|v-u|^{2H_0-2}$ with $c_H$ a constant depending only on $H$ and whose value can change from one line to another,  we have
\begin{eqnarray*}
&&\big(g_S(s,\cdot)\otimes_r g_T(t,\cdot)\big)(x_1,\ldots,x_{2q-2r})\\
&=&c_{H}\int_{0}^s \int_0^t  |v-u|^{(2H_0-2)r}
e^{-S(s-u)}e^{-T(t-v)}\prod_{j=1}^{q-r} (u-x_j)_+^{H_0-\frac32}\prod_{j=q-r+1}^{2q-2r}(v-y_j)_+^{H_0-\frac32}dudv.
\end{eqnarray*}
Now, let $\sigma,\gamma$ be two permutations of $\mathfrak{S}_{2q-2r}$, and
write $g_S(s,\cdot)\otimes_{\sigma,r} g_T(t,\cdot)$ to indicate
the function 
\[(x_1,\ldots,x_{2q-2r})\mapsto 
\big(g_S(s,\cdot)\otimes_r g_T(t,\cdot)\big)(x_{\sigma(1)},\ldots,x_{\sigma(2q-2r)}).
\]
We can write, for some integers $a_1,\ldots,a_4$ satisfying $a_1+a_2=a_3+a_4=q-r$ (and whose exact value is useless in what follows),
\begin{eqnarray*}
&&\big\langle
g_S(s,\cdot)\otimes_{\sigma,r} g_T(t,\cdot),
g_S(s,\cdot)\otimes_{\gamma,r} g_T(t,\cdot)
\big\rangle
\\
&=&c_{H}\int_{0}^s \int_0^t
\int_{0}^s \int_0^t 
 |v-u|^{(2H_0-2)r}
 |z-w|^{(2H_0-2)r}
 |u-w|^{(2H_0-2)a_1}
\\
 &&
\hskip3cm \times 
 |u-z|^{(2H_0-2)a_2}
 |v-w|^{(2H_0-2)a_3}
 |u-z|^{(2H_0-2)a_4}  
  \\
 &&\hskip3cm \times
e^{-S(s-u)}e^{-T(t-v)}
e^{-S(s-w)}e^{-T(t-z)}
dudvdwdz.
\end{eqnarray*}
We deduce that
\begin{eqnarray*}
&&(ST)^{\frac{2}{q}(1-H)+2H}
 \int_{[0,1]^2}\big\langle
g_S(s,\cdot)\otimes_{\sigma,r} g_T(t,\cdot),
g_S(s,\cdot)\otimes_{\gamma,r} g_T(t,\cdot)
\big\rangle dsdt\\
&=&c_H(ST)^{\frac{2}{q}(1-H)+2H}
 \int_{[0,1]^2}\left(
\int_{0}^s \int_0^t
\int_{0}^s \int_0^t 
 |v-u|^{(2H_0-2)r}
 |z-w|^{(2H_0-2)r}
\right.
\\
 &&
\hskip3cm \times 
 |u-w|^{(2H_0-2)a_1}
  |u-z|^{(2H_0-2)a_2}
 |v-w|^{(2H_0-2)a_3}
 |v-z|^{(2H_0-2)a_4}  
  \\
 &&\left.
 \hskip3cm \times
e^{-S(s-u)}e^{-T(t-v)}
e^{-S(s-w)}e^{-T(t-z)}
dudvdwdz\right)
dsdt\\
&=&c_H(ST)^{\frac{2}{q}(1-H)+2H}
 \int_{[0,1]^2}\left(
\int_{0}^s \int_0^t
\int_{0}^s \int_0^t 
 |v-u-t+s|^{(2H_0-2)r}
 |z-w+t-s|^{(2H_0-2)r}
\right.
\\
 &&
\hskip3cm \times 
 |u-w|^{(2H_0-2)a_1}
  |u-z+t-s|^{(2H_0-2)a_2}
 |v-w-t+s|^{(2H_0-2)a_3}
  \\
 &&\left.
 \hskip3cm \times
 |v-z|^{(2H_0-2)a_4}  
 e^{-Su}e^{-Tv}
e^{-Sw}e^{-Tz}
dudvdwdz\right)
dsdt\\
&=&c_HS^{\frac{2}q(1-H)(1+a_1-q)}
T^{\frac{2}q(1-H)(1+a_4-q)}\\
&&\hskip1.5cm\times
 \int_{[0,1]^2}\left(
\int_{0}^{Ss} \int_0^{Tt}
\int_{0}^{Ss} \int_0^{Tt} 
\left|\frac{v}T-\frac{u}{S}-t+s\right|^{(2H_0-2)r}
\left|\frac{z}T-\frac{w}S+t-s\right|^{(2H_0-2)r}
\right.
\\
 &&
\hskip3cm \times 
 |u-w|^{(2H_0-2)a_1}
 \left| \frac{u}S-\frac{z}T+t-s\right|^{(2H_0-2)a_2}
 \left| \frac{v}T-\frac{w}S-t+s\right|^{(2H_0-2)a_3}
  \\
 &&\left.
 \hskip3cm \times
 |v-z|^{(2H_0-2)a_4}  
 e^{-u}e^{-v}
e^{-w}e^{-z}
dudvdwdz\right)
dsdt.
\end{eqnarray*}

It follows that
\[
\lim_{S,T\to\infty} (ST)^{\frac{2}{q}(1-H)+2H}
 \int_{[0,1]^2}\big\langle
g_S(s,\cdot)\otimes_{\sigma,r} g_T(t,\cdot),
g_S(s,\cdot)\otimes_{\gamma,r} g_T(t,\cdot)
\big\rangle dsdt
\]
exists whatever $r$ and $a_1,\ldots,a_4$ such that $a_1+a_2=a_3+a_4=q-r$. Note that this limit is always zero, except when $r=1$, $a_1=a_4=q-1$ and $a_2=a_3=0$, in which case it is given by
\[
c_H \int_{[0,1]^2}
\left|t-s\right|^{4H_0-4}dtds\times 
\left(\int_{\R_+^2}
 |u-w|^{(2H_0-2)(q-1)}e^{-(u+w)}dudw\right)^2<\infty.
\]
Since
\[
g_S(s,\cdot)\widetilde{\otimes}_{r} g_T(t,\cdot) =
\frac{1}{(2q-2r)!}\sum_{\sigma\in\mathfrak{S}_{2q-2r}}
g_S(s,\cdot)\otimes_{\sigma,r} g_T(t,\cdot)
\]
the existence of the two limits (\ref{limit1})-(\ref{limit2}) follow, implying in turn the existence of $G_\infty$.
\medskip

Now, let us check the claim about the distribution of $G_\infty$.
Let $\widetilde{Y}_t=U_1(t)$, that is,
$\widetilde{Y}_t = \int_0^t e^{-(t-u)}dZ^{q, H}_u$, $t\geq 0$.
 By a scaling argument,
it is straightforward to check that
$(\widetilde{Y}_{tT})_{t\geq 0} \overset{\rm law}{=} T^H (U_T(t))_{t\geq 0}$
for any fixed $T>0$. As a result,
\[
T^{\frac{2}{q}(1-H)-1}\int_0^T (\widetilde{Y}_t^2-\E[\widetilde{Y}_t^2])dt = 
T^{\frac{2}{q}(1-H)}\int_0^1  (\widetilde{Y}_{tT}^2-\E[\widetilde{Y}_{tT}^2])dt 
\overset{\rm law}{=} G_T.
\]
Using (\ref{diu}), we deduce that $G_T/B_{H,q}$ converges in law
to the Rosenblatt distribution of parameter $1-\frac{2}{q}(1-H)$, hence the claim.\qed

\section{Proof of the consistency part in Theorem \ref{main}}

The consistency part of Theorem \ref{main} is directly obtained as a consequence of the following two propositions.

\begin{prop}
Let $X$ be given by (\ref{fracV})-(\ref{fracV2}) with $a>0$, $b\in\R$, $q\geq 1$ and $H\in(\frac12,1)$.
As $T\to\infty$, one has
\begin{equation}\label{lm41}
\frac{1}{T}\int_0^T X_tdt \to b\quad \mbox{a.s.}
\end{equation}
\end{prop}
{\it Proof}.
We use (\ref{fracV2}) to write
\[
\frac{1}{T}\int_0^T X_t dt = \frac{b}{T}\int_0^T (1-e^{-at}) dt + \frac{1}{T}\int_0^T Y_t dt.
\]
Since it is straightforward that $\frac{b}{T}\int_0^T (1-e^{-at}) dt\to b$, we are left to show that $\frac{1}{T}\int_0^T Y_t dt\to 0$ almost surely.

By (\ref{ivan2}), one can write, for any integer $n\geq 1$,
\[
\E\left[\left(\frac{1}{n}\int_0^n Y_t dt\right)^2\right]\leq \frac{2}{a^2n^2}\big(\E[(Z^{q, H}_n)^2]+\E[Y_n^2]\big) =O(n^{2H-2}),
\]
where the last equality comes from the $H$-selfsimilarity property of $Z^{q, H}$ as well as (\ref{O1}). 
Since $\frac{1}{T}\int_0^T Y_t dt$ belongs to the $q$th Wiener chaos, it enjoys the 
hypercontractivity property (\ref{eq:hypercontractivity1}). As a result, 
for all $p >\frac{1}{1-H}$ and $\lambda > 0$,
\begin{align*}
\sum_{n=1}^\infty \Prob\bigg(\Big|\frac{1}{n}\int_0^n Y_t dt\Big| > \lambda \bigg) 
&\leq \frac{1}{\lambda^p}\sum_{n=1}^\infty \E\bigg[\Big|\frac{1}{n}\int_0^n Y_t dt\Big|^p\bigg]
& \leq \frac{{\rm cst}(p)}{\lambda^p}\sum_{n=1}^\infty 
\E\bigg[\Big(\frac{1}{n}\int_0^nY_t dt\Big)^2\bigg]^{p/2}\\
& \leq   \frac{{\rm cst}(p)}{\lambda^p}\sum_{n=1}^\infty n^{-(1-H)p}<\infty.
\end{align*}
We deduce from the Borel-Cantelli lemma that
$
\frac{1}{n}\int_0^nY_t dt \rightarrow 0
$
almost surely as $n \to \infty$. 

Finally, fix $T > 0$ and let $n=\lfloor T\rfloor$ be its integer part. We can write
\begin{equation}\label{bla}
\frac{1}{T}\int_0^TY_t dt = \frac{1}{n}\int_0^nY_t dt + \frac{1}{T}\int_n^TY_t dt + \bigg(\frac{1}{T} - \frac{1}{n}\bigg)\int_0^nY_t dt.
\end{equation}
We have just proved above that $\frac{1}{n}\int_0^nY_t dt$ tends to zero almost surely as $n \to \infty$. We now consider the second and third terms in (\ref{bla}). 
We have, almost surely as $T\to\infty$,
 \[
 \bigg|\bigg(\frac{1}{T} - \frac{1}{n}\bigg)\int_0^nY_t dt \bigg| = \bigg(1- \frac{n}{T}\bigg) \bigg|\frac{1}{n}\int_0^nY_t dt \bigg| \leq \bigg|\frac{1}{n}\int_0^nY_t dt \bigg| \rightarrow 0,
 \]
and
\[
\bigg|\frac{1}{T}\int_n^TY_t dt\bigg| \leq \frac{1}{n} \int_n^{n+1}|Y_t|dt.
\]

To conclude, it remains to prove that $\frac{1}{n} \int_n^{n+1}|Y_t|dt \to 0$ almost surely as $n \to \infty$. Using (\ref{O1}) we have, for all fixed $\lambda > 0$, 
\begin{align*}
\mathbb{P} \bigg\{ \frac{1}{n} \int_n^{n+1}|Y_t|dt > \lambda \bigg\} & \leq \frac{1}{\lambda^2} \E\bigg[\bigg( \frac{1}{n} \int_n^{n+1}|Y_t|dt\bigg)^2\bigg] \\ 
& \leq \frac{1}{\lambda^2n^2}  \int_n^{n+1} \int_n^{n+1}  \sqrt{\E[Y_s^2]}\sqrt{\E[Y_t^2]}dsdt =O(n^{-2}).
\end{align*}
Hence, as $n \to \infty$, the Borel-Cantelli lemma applies and implies that $\frac{1}{n} \int_n^{n+1}|Y_t|dt$ goes to zero almost surely. 
This completes the proof of (\ref{lm41}).\qed

\begin{prop}
Let $X$ be given by (\ref{fracV})-(\ref{fracV2}) with $a>0$, $b\in\R, q\geq 1$ and $H \in (\frac12, 1)$. As $T\to\infty$, one has
\begin{equation}\label{lm42}
\frac{1}{T}\int_0^T X_t^2dt \to b^2+a^{-2H}H\Gamma(2H)\quad \mbox{a.s.}
\end{equation}
\end{prop}
\noindent
{\it Proof}. 
We first use (\ref{fracV2}) to write
\[
\frac{1}{T}\int_0^T X_t^2dt = \frac{1}{T}\int_0^T h(t)^2dt + \frac{2}{T}\int_0^T h(t)Y_tdt+\frac{1}{T}\int_0^T Y_t^2dt.
\]
We now study separately the three terms in the previous decomposition.
More precisely we will prove that,  as $T\to\infty$,
\begin{eqnarray}
\frac{1}{T}\int_0^T h(t)^2dt&\to& b^2,\label{dollar1}\\ 
\frac{1}{T}\int_0^T h(t)Y_tdt&\to& 0\quad\mbox{a.s.}\label{dollar2}\\ 
\frac{1}{T}\int_0^T Y_t^2dt&\to& a^{-2H}H\Gamma(2H)\quad\mbox{a.s.}\label{dollar3},
\end{eqnarray}
from which (\ref{lm42}) follows immediately.\\

\underline{First term}.
By Lebesgue dominated convergence, one has
\[
\frac{1}{T}\int_0^T h(t)^2dt = \int_0^1 h(Tt)^2dt =  b^2 \int_0^1 (1- e^{-aTt})^2 dt \rightarrow b^2,
\]
that is, (\ref{dollar1}) holds.

\bigskip

\underline{Second term}. 
First, we claim that
\begin{equation}\label{conv}
T^{-H}\int_0^T h(t)Y_tdt \overset{\rm law}{\to} \frac{b}{a} Z^{q, H}_1.
\end{equation}
Indeed, let us decompose:
\begin{eqnarray*}
\int_0^T h(t)Y_tdt=b\int_0^T (1-e^{-a t})Y_tdt
=b\int_0^T Y_tdt - b\int_0^T e^{-a t}\,Y_tdt.
\end{eqnarray*}
Using (\ref{O1}) in the last line, we can write
\begin{eqnarray}
&&\int_0^T e^{-at}\,Y_tdt = \int_0^T e^{-at}\left(\int_0^t e^{-a(t-s)}dZ^{q, H}_s\right)dt \notag\\
&=& \int_0^T \left(\int_s^T e^{-a(2t-s)}dt\right)dZ^{q, H}_s 
= \frac{1}{2a}\int_0^T (e^{-a(2T- s)}-e^{-as})dZ^{q, H}_s \notag\\
&=&\frac{1}{2a}\left( e^{-a T}Y_T - \int_0^T e^{-as}dZ^{q, H}_s\right)\notag\\
&\to&-\frac{1}{2a}\int_0^\infty e^{-a s}dZ^{q, H}_s\quad\mbox{in $L^2(\Omega)$ as $T\to\infty$.}\label{blabla}
\end{eqnarray}
The announced convergence (\ref{conv}) is a consequence of (\ref{ivan2}), (\ref{blabla}) and the selfsimilarity of $Z^{q, H}$.
Now, relying on the Borel-Cantelli lemma and the fact that $\int_0^T h(t)Y_tdt$ enjoys the
hypercontractivity property, it is not difficult to deduce from (\ref{conv}) that
(\ref{dollar2}) holds.\\
 
\underline{Third term}. Firstly, let us write, as $T\to\infty$,
\begin{align}\label{hou1}
\frac{1}{T}\int_0^T \E[Y_t^2]dt& = H(2H-1) \frac{1}{T}\int_0^Tdt \int_0^t\int_0^t dudv e^{-a u}e^{-a v}|u-v|^{2H-2} \nonumber \\
& = H(2H-1) \int_0^1 dt \int_0^{Tt}\int_0^{Tt} dudv e^{-a u}e^{-a v}|u-v|^{2H-2} \nonumber\\
&\longrightarrow  H(2H-1)  \int_0^\infty \int_0^\infty dudv e^{-a u}e^{-a v}|u-v|^{2H-2} \nonumber\\
&\quad\quad=a^{-2H}H\Gamma(2H). 
\end{align}
To conclude the proof of (\ref{dollar3}), we are thus left to show that :
\begin{equation}\label{eq:hou2}
\frac{1}{T}\int_0^T (Y_t^2 - \E[Y_t^2])dt \rightarrow 0 \text{ a.s.}
\end{equation}

First, we claim that, as $n \in \mathbb{N}^*$ goes to infinity,
\begin{equation}\label{eq:as1}
G_n:=\frac{1}{n}\int_0^n (Y_t^2 - \E[Y_t^2])dt \to 0 \text{ a.s.}
\end{equation}

Indeed, for all fixed $\lambda > 0$ and $p \geq 1$ we have, by the hypercontractivity property (\ref{eq:hypercontractivity2}) for $G_n$ belonging to a finite sum of Wiener chaoses,
\begin{align*}
\mathbb{P}\{ |G_n| > \lambda \}&\leq \frac{1}{\lambda^p} \E[ |G_n|^p ] \leq  \frac{\text{cst}(p)}{\lambda^p} \E[ G_n^2 ]^{p/2}.
\end{align*}

If ($q \geq 1$ and $H>\frac34$) or $q\geq 2$, combining (\ref{diu}) with, e.g., \cite[Lemma 2.4]{bNourdinPoly} leads to
\begin{equation}\label{azerty}
\sup_{T >0}  
\E\bigg[ 
\bigg(
T^{\frac2q(1-H)-1}\int_0^T (Y_t^2 - \E[Y_t^2])dt 
\bigg)^2 
\bigg] < \infty
\end{equation}
(note that one could also prove (\ref{azerty}) directly),
implying in turn that 
$
\mathbb{P}\{ |G_n| > \lambda \}
=O(n^{-\frac{2p}{q}(1-H)});$
choosing $p$ so that $\frac{2p}{q}(1-H) > 1$ leads to 
$ \sum_{n=1}^\infty  \mathbb{P}\{ |G_n| > \lambda \} < \infty$,
and so our claim (\ref{eq:as1}) follows from Borel-Cantelli lemma.

If $q=1$ and $H<\frac34$, the same reasoning (but using this time (\ref{ivan}) instead of (\ref{diu})) leads exactly to the same conclusion (\ref{eq:as1}).

Now, fix $T > 0$ and consider its integer part $n=\lfloor T\rfloor$. One has
\begin{equation}\label{blablabla}
G_T = G_n + \frac{1}{T}\int_n^T(Y_t^2 - \E[Y_t^2])dt + \bigg(\frac{1}{T} - \frac{1}{n}\bigg)\int_0^n(Y_t^2 - \E[Y_t^2])dt.
\end{equation}
We have just proved above that $G_n$ tends to zero almost surely as $n \to \infty$. We now consider the third term in (\ref{blablabla}). We have, using (\ref{eq:as1}):
\begin{align*}
\bigg|\frac{1}{T}- \frac{1}{n}\bigg|\bigg|\int_0^n(Y_t^2 - \E[Y_t^2])dt\bigg|&= \bigg(1-\frac{n}{T}\bigg)\bigg|\frac{1}{n}\int_0^n(Y_t^2 - \E[Y_t^2])dt\bigg|
\leq  |G_n| \rightarrow 0 \quad \text{a.s}.
\end{align*}
Finally, as far as the second term in (\ref{blablabla}) is concerned, we have
\[
\bigg|\frac{1}{T}\int_n^T(Y_t^2- \E[Y_t^2])dt\bigg|  \leq \frac{1}{n}\int_n^{n+1}|Y_t^2 - \E[Y_t^2]|dt.
\]
To conclude, it thus remains to prove that, as $n \to \infty$,
\begin{equation}\label{eq:10}
F_n:=\frac{1}{n}\int_n^{n+1}|Y_t^2 - \E[Y_t^2]|dt \to 0 \qquad \text{a.s}.
\end{equation}
By hypercontractivity and (\ref{yt2}), one can write
\begin{align*}
 \text{Var}(Y_t^2)  \leq \text{cst}(q)(\E[Y_t^2])^2 \leq \text{cst}(q)  a^{-4H}H^2\Gamma(2H)^2.
\end{align*}
Thus, $\sup_{t} \text{Var}(Y_t^2) < \infty,$ and it follows that
\[
\E[F_n^2] = \frac{1}{n^2}\int_n^{n+1}\int_n^{n+1}\E\big[ \big|Y_t^2 - \E[Y_t^2]\big|\big|Y_s^2 - \E[Y_s^2]\big| \big]dsdt =O(n^{-2}).
\]
Hence $ \sum_{n=1}^\infty  \mathbb{P}\{ |F_n| > \lambda \} \leq \sum_{n=1}^\infty \frac{1}{\lambda^2}E[F_n^2] < \infty$ for all $\lambda > 0$,
and Borel-Cantelli lemma leads to (\ref{eq:10}) and concludes the proof of (\ref{dollar3}).\qed

\section{Proof of the fluctuation part in Theorem \ref{main}}

We now turn to
the proof of the part of Theorem \ref{main} related to fluctuations.
We start with the fluctuations of $\widehat{b}_T$, which are easier compared to
$\widehat{a}_T$.\\

\underline{Fluctuations of $\widehat{b}_T$}. Using first (\ref{decompoX}) and then (\ref{O1}) and (\ref{ivan2}),
we can write
\begin{equation}\label{bt}
T^{1-H}\big\{\widehat{b}_T - b\} =T^{1-H}\left\{\frac1T\int_0^T Y_tdt - \frac{b}T\int_0^T e^{-at}dt \right\} 
=\frac{Z^{q, H}_T}{aT^H} +O(T^{-H}),
\end{equation}
which will be enough to conclude, see the end of the present section.

\medskip

\underline{Fluctuations of $\widehat{a}_T$}. 
As a preliminary step, we first concentrate on the asymptotic behavior, as $T\to\infty$, of the random quantity 
\[
\ell_T:=\alpha_T -a^{-2H}H\Gamma(2H),
\] 
where $\alpha_T$ is given by (\ref{alpha_T}).
Since $X_t=h(t)+Y_t$, see (\ref{decompoX}),  we have
\begin{eqnarray*}
\ell_T 
=A_T+B_T+2C_T+D_T-E_T^2-2E_TF_T-F_T^2,
\end{eqnarray*}
where
\begin{eqnarray*}
A_T&=&\frac{1}{T}\int_0^T (Y_t^2-\E[Y_t^2])dt,\quad
B_T=\frac{1}{T}\int_0^T \E[Y_t^2]dt - a^{-2H}H\Gamma(2H)\\
C_T&=&\frac{1}{T}\int_0^T Y_t h(t)dt,\quad D_T=\frac{1}{T}\int_0^T h^2(t)dt,\quad
E_T=\frac{1}{T}\int_0^T Y_tdt,\quad F_T=\frac{1}{T}\int_0^T h(t)dt.
\end{eqnarray*}
We now treat each of these terms separately.
\medskip

\noindent
{\it Term $B_T$}. Recall from (\ref{yt2}) that, as $T \to \infty$,
\begin{eqnarray*}
\E[Y_T^2]&=&H(2H-1)\int_{[0,T]^2}e^{-a(u+v)}|u-v|^{2H-2}dudv\\
&\to& H(2H-1) \int_{[0,\infty)^2}e^{-a(u+v)}|u-v|^{2H-2}dudv = a^{-2H}H\Gamma(2H).
\end{eqnarray*} 
As a result,
\begin{eqnarray*}
|B_T|&=&\left|\frac{1}{T}\int_0^T \E[Y_t^2]dt - a^{-2H}H\Gamma(2H)\right|\\
&\leq&\frac{H(2H-1)}{T}\int_0^T dt \int_{[0,\infty)^2\setminus[0,t]^2} dudv\,e^{-a (u+v)}|u-v|^{2H-2}\\
&\leq&\frac{2H(2H-1)}{T}\int_0^T dt\int_t^\infty dv\,e^{-a v}\int_0^\infty du\,e^{-a u}{\bf 1}_{\{v\geq u\}}(v-u)^{2H-2}\\
&\leq&\frac{2H(2H-1)}{T}\int_0^\infty dt\int_t^\infty dv\,e^{-a v}\int_0^v du\,u^{2H-2}\\
&=&\frac{2H}{T}\int_0^\infty dt\int_t^\infty dv\,e^{-a v}v^{2H-1}
=\frac{2H}{T}\int_0^\infty e^{-a v}v^{2H}dv=O(\frac1T).
\end{eqnarray*} 

\medskip
\noindent
{\it Term $C_T$}. We can write
\[
C_T=\frac{1}{T}\int_0^T Y_t h(t)dt =\frac{b}T\int_0^T (1-e^{-a t})Y_tdt
=b\left(\frac1T\int_0^T Y_tdt - \frac1T\int_0^T e^{-a t}Y_tdt\right).
\]
But
\begin{eqnarray*}
\frac1T\int_0^T e^{-a t}Y_tdt&=&
\frac1T\int_0^T e^{-a t}\left(\int_0^t e^{-a(t-s)}dZ^{q, H}_s\right)dt\\
&=&\frac1T\int_0^T e^{a s}\left(
\int_s^T e^{-2a t}dt
\right)dZ^{q, H}_s=\frac1{2a T}\left(\int_0^T e^{-a s}dZ^{q, H}_s - e^{-a T}Y_T\right).
\end{eqnarray*}
Using  (\ref{O1}) and $\int_0^T e^{-a s}dZ^{q, H}_s\to \int_0^\infty e^{-a s}dZ^{q, H}_s$ in $L^2(\Omega)$, we deduce that
\[
C_T = b\,E_T+O(\frac1T).
\]
\medskip
\noindent
{\it Term $D_T$}. It is straightforward to check that 
\[
D_T=\frac{1}{T}\int_0^T h^2(t)dt = \frac{b^2}{T}\int_0^T (1-e^{-a t})^2dt=b^2+O(\frac1T).
\]

\medskip
\noindent
{\it Term $E_T$}. Thanks to (\ref{O1}) and (\ref{ivan2}), we have
\[
E_T = \frac1T\int_0^T Y_Tdt =  \frac{1}{a\,T}(Z^{q, H}_T-Y_T)=
\frac{Z^{q, H}_T}{a\,T} + O(\frac1T).
\]
Since $Z^{q, H}_T\overset{\rm law}{=}T^HZ^{q, H}_1$ by selfsimilarity, we deduce
\[
E_T^2
=\frac{(Z^{q, H}_T)^2}{a^2T^2}+ O(T^{H-2}). 
\]
\medskip
\noindent
{\it Term $F_T$}. Similarly to $D_T$, it is straightforward to check that
\[
F_T=\frac{1}{T}\int_0^T h(t)dt = \frac{b}{ T}\int_0^T (1-e^{-a t})dt=b+O(\frac1T).
\]

\medskip
\noindent
Combining everything together, we eventually obtain that
\begin{equation}\label{comb}
\ell_T = A_T-\frac{(Z^{q, H}_T)^2}{a^2T^2}+O(T^{-1}).
\end{equation}

\medskip

\underline{Fluctuations of $(\widehat{a}_T,\widehat{b}_T)$}. 
$\bullet$ We consider first the case $(q=1 \text{ and } H > \frac34)$ or $(q \geq 2)$. Since $Z^{q, H}$ satisfies the scaling property, we can write
\[
\Big( T^{\frac2q(1-H)}A_T, 
T^{-H}Z^{q, H}_T \Big)
\overset{\rm law}{=}
\left(
a^{-\frac2q(1-H)-2H}\,G_{aT},
Z^{q, H}_1
\right),
\]
and we deduce from the $L^2$-convergence of $G_T$ (see Proposition \ref{finfini})
that
\begin{equation}\label{eq:euro}
\Big( T^{\frac2q(1-H)}A_T, 
T^{-H}Z^{q, H}_T \Big)
\overset{\rm law}{\to}
\left(a^{-\frac2q(1-H)-2H}\,G_\infty,
Z^{q, H}_1\right).
\end{equation}
On the other hand, a Taylor expansion yields
\begin{eqnarray*}
T^{\frac2q(1-H)}\{\widehat{a}_T-a\}&=&T^{\frac2q(1-H)}\,a\left[\left(1+\frac{a^{2H}\,
\ell_T }{H\Gamma(2H)}\right)^{-\frac1{2H}}-1\right]\\
&=& -\frac{a^{1+2H}}{2H^2\Gamma(2H)}\,\left(T^{\frac2q(1-H)}A_T
-T^{\frac2q(1-H)-2}\frac{(Z^{q, H}_T)^2}{a^2}\right)+o(1),
\end{eqnarray*}
implying in turn
\begin{eqnarray*}
&&\left(T^{\frac2q(1-H)}\{\widehat{a}_T-a\},T^{1-H}\big\{\widehat{b}_T - b\} \right)\\
&=&
\left(-\frac{a^{1+2H}}{2H^2\Gamma(2H)}\,\left[T^{\frac2q(1-H)}A_T
-T^{2(\frac1q-1)(1-H)}\frac{(T^{-H}Z^{q, H}_T)^2}{a^2}\right],\frac{T^{-H}Z^{q, H}_T}{a}\right)+o(1),
\end{eqnarray*}
so that, from (\ref{eq:euro}),
\begin{eqnarray*}
&&\left(T^{\frac2q(1-H)}\{\widehat{a}_T-a\},T^{1-H}\big\{\widehat{b}_T - b\} \right)\\
& \overset{\rm law}{\to}&
\left\{
\begin{array}{lll}
\left(-\frac{a^{1-\frac2q(1-H)}}{2H^2\Gamma(2H)}\,G_\infty,\frac{Z^{q, H}_1}{a}\right)&\mbox{ if $q\geq 2$}\\
\quad\\
\left(-\frac{a^{2H-1}}{2H^2\Gamma(2H)}\,(G_\infty- (B^H_1)^2),\frac{B^H_1}{a}\right)&\mbox{ if $q=1$ and $H>\frac34$}
\end{array}
\right.,
\end{eqnarray*}
as claimed.

$\bullet$ Assume now that $q=1$ and $H< \frac34$, and let us write $B^H$ instead of $Z^{1,H}$ for simplicity. 
We deduce from (\ref{comb}) and $T^{-2}(B^{H}_T)^2\overset{\rm law}{=} T^{2H-2}(B^{ H}_1)^2$ that 
$
\sqrt{T}\ell_T = \sqrt{T}A_T+o(1),
$
so that
\begin{eqnarray*}
\left(\sqrt{T}\{\widehat{a}_T-a\},T^{1-H}\big\{\widehat{b}_T - b\} \right)
=
\left(-\frac{a^{1+2H}}{2H^2\Gamma(2H)}\,\sqrt{T}A_T,\frac{T^{-H}B^H_T}{a}\right)+o(1),
\end{eqnarray*}
implying in turn by (\ref{ivanbis}) that
\begin{eqnarray*}
\left(\sqrt{T}\{\widehat{a}_T-a\},T^{1-H}\big\{\widehat{b}_T - b\} \right)\to
\left(-\frac{a^{1+4H}\sigma_H}{2H^2\Gamma(2H)}\, N,\frac{N'}{a}\right),
\end{eqnarray*}
where $N,N'\sim N(0,1)$ are independent, as claimed.

$\bullet$ Finally, we consider the case $q=1$ and $H= \frac34$.
We deduce again from (\ref{comb}) and $T^{-2}(B^{3/4}_T)^2\overset{\rm law}{=} T^{-\frac12}(B^{3/4}_1)^2$ that 
\[
\sqrt{\frac{T}{\log T}}\ell_T = \sqrt{\frac{T}{\log T}}A_T+o(1),
\]
so that, using (\ref{nunubis}),
\begin{eqnarray*}
\left(\sqrt{\frac{T}{\log T}}\{\widehat{a}_T-a\},T^{\frac14}\big\{\widehat{b}_T - b\} \right)
&=&
\left(-\frac{16a^{\frac52}}{9\sqrt{\pi}}\,\sqrt{\frac{T}{\log T}}A_T,\frac{T^{-\frac34}B^H_T}{a}\right)+o(1)\\
&\to&
\left(\frac{3}4\sqrt{\frac{a}\pi}\, N,\frac{N'}{a}\right),
\end{eqnarray*}
where $N,N'\sim N(0,1)$ are independent.

\qed

{\bf Acknowledgment}. 
We gratefully acknowledge an anonymous referee  for a very constructive report, which led to a significant improvement of the paper, especially its link with the financial literature.

\end{document}